\documentclass[12pt]{article}
\usepackage{amssymb,amsmath}

\usepackage{graphicx}
\usepackage[arrow,matrix,curve]{xy}

\usepackage{amscd}
\usepackage{amsfonts}
\usepackage{amsthm}
\usepackage{mathrsfs}

\DeclareMathAlphabet{\mathpzc}{OT1}{pzc}{m}{it}

\newcommand{\e}{{\rm e}}
\newcommand{\ii}{{\rm i}}
\newcommand{\eqn}[1]{(\ref{#1})}

\newcommand{\real}{{\mathbb R}} 
\newcommand{\bea}{\begin{eqnarray}}
\newcommand{\eea}{\end{eqnarray}}
\def\beqa{\begin{eqnarray}}
\def\eeqa{\end{eqnarray}}
\def\nn{\nonumber}
\newcommand{\del}{\partial}

\newcommand{\eq}{\begin{equation}}
\newcommand{\eqa}{\begin{eqnarray}}
\newcommand{\en}{\end{equation}}
\newcommand{\ena}{\end{eqnarray}}

\def\sk{\vskip .4cm}
\def\noi{\noindent}
\def\nn{\nonumber}
\def\epsi{{\varepsilon}}
\def\f{{\rm{f}\,}}

\newcommand{\om}{\omega}
\newcommand{\Om}{\Omega}
\newcommand{\al}{\alpha}
\newcommand{\la}{\lambda}

\newcommand{\Ga}{\Gamma}

\newcommand{\RR}{{\mathcal R}}
\newcommand{\oR}{{\bar{\R}}}
\newcommand{\st}{{\star}}

\newcommand{\UU}{{U}\Xi}
\newcommand{\FF}{\mathcal F}
\newcommand{\varepsi}{\varepsilon}
\newcommand{\A}{{\mathcal  A}}
\newcommand{\AAs}{{A_\st }}
\newcommand{\UUs}{{\UU_\st }}
\newcommand{\Xis}{{\Xi_\st }}
\newcommand{\Oms}{\Omega_\st}
\newcommand{\La}{\Lambda }

\newcommand{\D}{\Delta}
\newcommand{\TT}{{\mathcal T}}

\newcommand{\re}{{\rangle}}
\newcommand{\OO}{\mathsf{\Omega}}
\newcommand{\TA}{\mathsf{\Theta}}

\newcommand{\ots}{\otimes_\st}
\newcommand{\res}{\re_\st}

\newcommand{\x}{\mathsf{x}}
\newcommand{\y}{\mathsf{y}}
\newcommand{\cc}{\mathbb{C}}

\newcommand{\ff}{{\sf f}}
\newcommand{\rr}{{\R}}

\def\ll{{\mathcal L}}
\def\D{\Delta}
\def\st {\star}
\def\dd{{\triangledown}}
\def\dds{\triangledown^\st}

\def\be{\beta}

\def\of{{\bar{{\rm{f}}\,}}}
\def\R{{R}}
\def\oR{{\bar{\R}}}
\def\cc{\mathbb{C}}

\def\TT{{\mathcal T}} 
\def\Xis{{\Xi_\st }}
\def\ots{\otimes_\st}
\def\rr{\mathsf{R}}
\def\tr{\mathsf{T}}
\def\XD{\mathsf{X}}
\def\Oms{\Omega_\st}
\def\Om{\Omega}

\begin{document}
\begin{titlepage}
\begin{flushright}

\baselineskip=12pt
DISTA-UPO/08\\
\hfill{ }\\
\end{flushright}

\begin{center}

\baselineskip=24pt

{\Large\bf Star Product Geometries}

\baselineskip=14pt

\vspace{1cm}

{\bf Paolo Aschieri}
\\[6mm]
{\it Centro Studi e Ricerche ``Enrico Fermi'' Compendio Viminale, 00184 Roma, Italy}\\
 and  {\it Dipartimento di Scienze e Tecnologie
 Avanzate, Universit\`{a} del
 Piemonte Orientale,}\\ and {\it INFN, Sezione di Torino, gruppo collegato di Alessandria }\\
aschieri@to.infn.it
\\[10mm]

\end{center}

\vskip 2 cm

\begin{abstract}
We consider noncommutative geometries obtained from a triangular
Drinfeld twist.
This allows to construct and study a wide class of noncommutative manifolds and their 
deformed Lie algebras of infinitesimal diffeomorphisms. 
This way 
symmetry principles can be implemented. We review two main examples \cite{GR1}-\cite{ALV}:  
a) general covariance in noncommutative spacetime. This leads to a noncommutative gravity theory. b) Symplectomorphims of the algebra of observables associated to a noncommutative configuration space. This leads to a geometric formulation of quantization on noncommutative 
spacetime, i.e., we establish a noncommutative correspondence principle from
$\st$-Poisson brackets to $\st$-commutators.
\vskip .25 cm 
New results concerning noncommutative gravity include the Cartan structural equations 
for the torsion and curvature tensors, and the associated Bianchi identities. Concerning 
scalar field theories the deformed algebra of classical and quantum observables has been understood in terms of a twist $\widehat{\mathscr{F}}$ within the algebra. 

\end{abstract}

\end{titlepage}

\section{Introduction}
An interesting and promising field of research
is the issue of spacetime structure in extremal energy regimes. 
There are evidences from General
Relativity, string theory and black hole physics which support the
hypothesis of a noncommutative structure.  The simplest and
probably most suggestive argument which points at a failure of the
classical spacetime picture at high energy scales comes from the
attempt of conjugating the principles of Quantum Mechanics with
those of General Relativity (see \cite{DoplicherFredenhagenRoberts},
and for a review \cite{doplicher}). If
one tries to locate an event with a spatial accuracy comparable
with the Planck length ($\sim 10^{-33}$ cm), spacetime uncertainty relations
necessarily emerge. In total analogy with Quantum Mechanics,
uncertainty relations are naturally implied by the presence of
noncommuting coordinates,
$[x^\mu,x^\nu]\not=0$.
Therefore,
below Planck length the usual description of spacetime as a
pseudo-Riemannian manifold locally modeled on Minkowski space is not
adequate anymore, and it has been proposed that it be described by
a \emph{Noncommutative Geometry} \cite{connes,landi,ticos}. This
line of thought has been pursued since the early days of Quantum
Mechanics \cite{heisenberg}.

In this context two relevant issues are the formulation of General 
Relativity and the quantization of field theories on noncommutative 
spacetime. We address these issues by developing a
differential geometry on noncommutative spacetime and phase space.
\sk
There are many noncommutative differential geometries (Connes type,
quantum groups like, matrix or fuzzy,...).  
We work in the deformation quantization context; noncommutativity
is obtained by introducing a formal $\st$-product on the algebra of
smooth functions. The most widely studied form of
noncommutativity is the Moyal-Weyl one (on spacetime $\real^4$),
$x^\mu\st x^\nu-x^\nu\st x^\mu=\ii\theta^{\mu\nu}$ 
with $\theta^{\mu\nu}$ a constant antisymmetric matrix. This
noncommutativity is obtained trought the $\star$-product
$
\left(f\star h\right)(x) = {\rm exp} \big(\frac{\ii}{2}
\theta^{\mu\nu}\frac{\del}{\del x^\mu} \frac{\del}{\del y^\nu}\big) 
f(x)h(y)\big|_{x=y}
~.
$

We focus on a class of star products obtained via a triangular 
Drinfeld twist \cite{drinfeld}. 
This is a quite wide class (the examples in Section 3 include  
quantization of symplectic and also of Poisson structures).
The algebra of functions of the
noncommutative  torus, of the noncommutative spheres 
\cite{Connes-Landi} and of further noncommutative manifolds 
(so-called isospectral deformations) considered in 
\cite{Connes-Landi}, and in \cite{Dubois-Violette}, 
\cite{Gayral:2005ad}, is associated to a 
$\st$-product structure obtained via a triangular Drinfeld twist
(see
\cite{Sitarz} and, for the four-sphere in \cite{Connes-Landi},
see \cite{Varilly}, \cite{Aschieri-Bonechi}).   
The star products we study are however not the most general ones, in 
particular they are a subclass of those associated with a 
quasitriangular structure \cite{drinfeld2}: on that noncommutative 
algebra of functions there is an action of the braid group, 
in the case we study there is an action of the permutation group.

It is remarkable how far in the program of formulating a noncommutative
differential geometry one can go using triangular Drinfeld twists. 
The study of this
class of $\st$-products geometries are first examples that can uncover 
some common features of a wider class of noncommutative geometries.

\sk
In Section 2 we introduce the twist $\FF$
and, starting from the principle that
every product, and in general every bilinear map, 
is consistently deformed by composing it with the appropriate realization of 
the twist $\mathcal{F}$, we review the construction of noncommutative 
differential geometry as in \cite{GR2,GR1,AschieriCorfu}.
A key point is that vectorfields have a  natural $\st$-action 
($\st$-Lie derivative) on the noncommutative algebras of functions  
and tensorfields, giving rise to the concept of deformed derivations.
These $\st$-derivations form a quantum Lie algebra ($\st$-Lie algebra). 
In this way we consider the $\st$-Lie algebra of vectorfields.
Next we can define covariant derivatives along vectorfields
because their deformed Leibniz rule is the same as for Lie derivatives 
along vectorfields and 
because covariant derivative and Lie derivative coincide on functions. 
New results is this section include the Cartan structural equations 
for the torsion and curvature tensors, and the associated Bianchi identities. 

In Section 3, following and developing the results of \cite{ALV}  we consider a noncommutative manifold 
$M$ (described by a 
twist $\FF$ and an associated $\st$-product) with an extra 
Poisson structure $\La$. 
If the noncommutative and the Poisson structures are compatible
(in particular if the Poisson structure $\La$ and the Poisson 
structure arising from the semiclassical limit of the $\st$-product 
are compatible) then the Poisson algebra of functions on $M$
can be deformed in a noncommutative $\st$-Poisson algebra.
Otherwise stated 
the algebra of observables becomes noncommutative and 
equipped with a compatible $\st$-Lie algebra structure. 
Correspondingly
the $\st$-Lie algebra of 
Hamiltonian vectorfields is a subalgebra of the $\st$-Lie algebra of 
vectorfields.

In Section 4 we apply this general construction to the 
infinite dimensional phase space of a scalar field theory. Here 
phase space has the usual Poisson structure $\{\Phi(x),\Pi(y)\}=\delta(x-y)$,
it also has a noncommutative structure ($\FF$ and its $\st$-product) 
induced from spacetime noncommutativity.
Usual $\hbar$-quantization of this phase space can also be deformed with 
the twist $\FF$ ($\st$-deformed).
We thus construct a $\st$-deformation of the algebra and Lie algebra of quantum 
observables of a scalar field. Starting from the usual canonical quantization 
map for field theories on commutative spacetime, $\Phi\stackrel{\hbar}
{\rightarrow}\hat\Phi$, we thus uniquely obtain  \cite{ALV} a quantization 
scheme for field theories on noncommutative spacetime, and show that it 
satisfies a correspondence principle between $\st$-Poisson brackets and 
$\st$-commutators.

\section{Deformation by twists}
In this section we describe the general setting 
used to introduce a star product via a twist.

Consider a Lie algebra $g$ over $\cc$ , and its associated universal enveloping algebra $Ug$.  We recall that the elements of $Ug$ are 
the complex numbers $\cc$ and sums 
of products of elements $t\in g$, where we identify 
$tt'-t't$ with the Lie algebra element $[t,t']$. 
$Ug$ is an associative algebra with unit. It is a Hopf algebra with coproduct 
$ \Delta:Ug\rightarrow Ug\otimes Ug~$, counit $\epsi : Ug\rightarrow \cc$
and antipode $S$ given on the generators as:
\begin{eqnarray}
 & & \D(t)=t \otimes 1 + 1 \otimes t~~~\D(1)=1\otimes 1  
  \label{copL}\\
 & & \epsi (t)=0~~~~~~~~~~~~~~~~~~~~\epsi (1)=1 \label{couL}\\
 & & S(t)=-t~~~~~~~~~~~~~~~~~S (1)=1 \label{coiL}
\end{eqnarray}
and extended to all $U(g)$ by requiring 
$\D$ and $\epsi$ to be linear and multiplicative
(e.g. $\D(tt'):=\D(t)\D(t')=tt'\otimes 1+t\otimes t' +t'\otimes t+ 1\otimes tt'$),
while $S$ is linear and antimultiplicative.
In the sequel we use Sweedler coproduct notation
\eq
\Delta(\xi)=\xi_1\otimes \xi_2
\en
where $\xi\in Ug$,  $\xi_1\otimes \xi_2\in Ug\otimes Ug$ and a sum
over $\xi_1$ and $\xi_2$  is understood.

We extend the notion of enveloping algebra to formal power series in 
$\la$ (we replace the field $\cc$ with the ring $\cc[[\la]]$)  and we 
correspondingly consider the Hopf algebra $(Ug[[\lambda]],\cdot,\Delta, S,\varepsi)$.
In the sequel for sake of brevity we will often denote 
$Ug[[\la]]$ by $Ug$.
\sk
A twist $\FF$ is an element 
$\FF\in Ug[[\la]]\otimes Ug[[\la]]$ that is invertible and that satisfies
\eq\label{propF1}
\FF_{12}(\Delta\otimes id)\FF=\FF_{23}(id\otimes \Delta)\FF\,,
\en
\eq\label{propF2}
(\varepsi\otimes id)\FF=1=(id\otimes \varepsi)\FF~,
\en
where $\FF_{12}=\FF\otimes 1$ and $\FF_{23}=1\otimes \FF$.

In our context we in addition require\footnote{ Actually it is
possible to show that (\ref{consF}) is a consequence  of
(\ref{propF1}), (\ref{propF2}) and of $\FF$ being at each order in
$\la$ a finite sum of finite products of Lie algebra elements} 
\eq
\label{consF}
\FF=1\otimes 1 + {\cal O}(\la)~.
\en
Property (\ref{propF1}) states that $\FF$ is a two cocycle, and it will turn out to be responsible for the associativity of the $\st$-products to be defined. Property 
(\ref{propF2}) is just a
normalization condition. From (\ref{consF}) it follows that 
${\cal F}$ can be formally inverted as a power series in $\lambda$. It also shows that the geometry we are going to construct has the nature of a deformation, i.e. in the $0$-th order in $\lambda$ we recover  the usual undeformed geometry.

\sk
 
We shall frequently use the notation (sum over $\al$ understood)
\eq\label{Fff}
\FF=\f^\al\otimes\f_\al~~~,~~~~\FF^{-1}=\of^\al\otimes\of_\al~.
\en
The elements $ \f^\alpha, \f_\alpha, \of^\al,\of_\al$ live in  $Ug$. 

In order to get familiar with this notation we  rewrite equation
(\ref{propF1}) and its inverse, 
\eq\label{ifpppop}
((\Delta\otimes id)\FF^{-1}) 
\FF^{-1}_{12} =((id \otimes \Delta)\FF^{-1})\FF^{-1}_{23}~,
\en 
as well as (\ref{propF2}) and (\ref{consF}) using the notation (\ref{Fff}), explicitly
\begin{eqnarray}
\f^\beta \f^\alpha_{_1}\otimes \f_\beta \f^\alpha_{_2}\otimes \f_\alpha &=& 
\f^\alpha\otimes \f^\beta \f_{\alpha_1}\otimes \f_\beta
\f_{\alpha_2} ,\label{2.21}\\
\label{ass}
\of_{_1}^\al\of^\be\otimes \of_{_2}^\al\of_\be\otimes \of_\al&=&
\of^\al\otimes {\of_{\al_1}}\of^\be\otimes {\of_{\al_2}}\of_\be~,\\
\varepsilon(\f^\alpha)\f_\alpha &=& 1 =  \f^\alpha\varepsilon(\f_\alpha) , \label{2.23}\\
{\cal F} = \f^\alpha\otimes \f_\alpha &=& 1\otimes 1 + {\cal O}(\lambda) .\label{2.24}
\end{eqnarray}

\sk
Consider now an algebra $A$ (over $\cc[[\la]$), 
and an action of the Lie algebra $g$ on $A$, $a\mapsto t(a)$ where $t\in g$ 
and $a\in A$. We require compatibility of this action with the product in $A$
i.e.  for any $t\in g$ we have a derivation of $A$, 
\eq
t(ab)=t(a)b+at(b)~.\label{compprim}
\en
The action of $g$ on $A$ induces an action of the universal enveloping 
algebra $Ug$ on $A$  (for example the element $tt'\in Ug$ has action 
$t(t'(a))$). We say that $A$ is a $Ug$-module algebra, i.e.,
the  algebra structure of the $Ug$-module $A$ 
is compatible with the $Ug$ action,
for all $\xi\in \UU$ and $ a,b\in A$,
\eq\xi(ab)=\mu\circ\Delta(\xi)(a\otimes b)=\xi_1(a)\xi_2(b)~~~,~~~~
\,\xi(1)=\varepsilon(\xi)1\,. \label{Ugmodulealg}
\en
(where $1$ is the unit in $A$). This property is equivalent to  
\eqn{compprim}. 
\sk
Given a twist $\FF\in Ug\otimes Ug$,
we can construct a deformed algebra $A_\st$.
The algebra $A_\st$ has the same vector space 
structure as $A$. 
The product in $A_\st $ is defined by
\eq
a\st  b=\mu\circ \FF^{-1}(a\otimes b)=\of^\al(a)\of_\al(b)~.
\en
In order to prove associativity of the new product we
use (\ref{ass}) and compute:
\eqa
(a\st b)\st c &=&\of^\al(\of^\be(a)\of_\be(b))\of_\al(c)=
(\of_{_1}^\al\of^\be)(a) (\of_{_2}^\al\of_\be)(b) \of_\al(c)\nn
=
\of^\al(a) ({\of_{\al_1}}\of^\be)(b)
({\of_{\al_2}}\of_\be)(c)\nn \\
&=&
\of^\al(a)\of_\al(\of^\be(b)\of_\be(c))=a\st (b\st c)\nn~.
\ena

We will see that $A_\st$ is a module algebra with respect to
a deformed Hopf algebra $Ug_\st$ or the isomorphic Hopf algebra $Ug^\FF$,
(cf. text after \eqn{befmod}, and  Note 1 at the end of this section).

\sk
\sk
\noi{\bf{Examples of twists}}\\ 
Consider the element in $Ug[[\la]]\otimes Ug[[\la]]$ given by
\eq
\FF=e^{-{i\over 2}\la\theta^{\mu\nu}t_\mu
\otimes t_\nu}
\en 
where the elements $\{t_\mu\}$ generate an abelian subalgebra of $g$, and
$\theta^{\mu\nu}$ is a constant matrix (usually antisymmetric and real
in order to be compatible with conjugation).
The inverse of $\FF$ is
$$\FF^{-1}=e^{{i\over 2}\la\theta^{\mu\nu}t_\mu
\otimes t_\nu~.}$$
Then we have
$(\Delta\otimes id)\FF=e^{-{i\over 2}\la\theta^{\mu\nu}
(t_\mu\otimes 1\otimes t_{\nu}+
1\otimes t_\mu
\otimes t_\nu)}$
so that property (\ref{propF1}) easily follows: 
$\FF_{12}(\Delta\otimes id)\FF=e^{-{i\over 2}\la\theta^{\mu\nu}(t_\mu
\otimes t_\nu\otimes 1+
t_\mu\otimes 1
\otimes t_\nu +
1\otimes t_\mu
\otimes t_\nu)}
=\FF_{23}(id\otimes \Delta)\FF\,.$
Property (\ref{propF2}) trivially holds.
\sk
Twists are not necessarily related to abelian Lie algebras. 
For example consider
the elements $H, E,A,B$, that satisfy the Lie algebra relations
\eqa
& & [H,E]=2E~,~~ [H,A]=\al A~,~~[H,B]=\be B~,~~~~~~~~\al+\be=2~,\nn\\[.1cm]
& & [A,B]=E~,~~~~[E,A]=0~,~~~\;~[E,B]=0~.   
\ena
Then the element
\eq
\FF=e^{{1\over 2}H\otimes_{} {\rm{ln}}(1+\la E)} \:  e^{\la A\otimes B{1\over 1+\la E}}
\en
is a twist and gives a well defined $\st$-product on the algebra of
functions on $M$.
These twists are known as extended Jordanian deformations
\cite{Kulish1}. Jordanian deformations 
\cite{GER, OGIEV} are obtained setting $A=B=0$ (and keeping the relation 
$[H,E]=2E$).

\subsection{$\st$-Noncommutative Manifolds}

\noi  {\it Definition 1. } A $\st$-noncommutative manifold is a quadruple 
$(M,g,\FF,\rho)$ where $M$ is a smooth manifold, $g$ is a Lie algebra 
with twist $\FF\in Ug[[\la]]\otimes Ug[[\la]]$ and $\rho :g\rightarrow \Xi$ is a homomorphism of $g$ in the Lie algebra $\Xi$ of vectorfields on $M$.
\sk
Usually we just consider the image $\rho(g)$ of $g$ in $\Xi$. Then the $\st$-noncommutative manifold is defined by the couple $(M,\FF)$ where 
$\FF\in U\Xi[[\la]]\otimes U\Xi[[\la]]$. 

We now use the twist to  deform the commutative geometry on a manifold $M$
(vectorfields, 1-forms, exterior algebra, tensor algebra, symmetry algebras, 
covariant derivatives etc.) 
into the twisted noncommutative one. 
The guiding principle is the observation that every time we have 
a bilinear map $$\mu\,: X\times Y\rightarrow Z~~~~~~~~~~~~~~~~~~$$
where $X,Y,Z$ are vectorspaces, 
and where there is an action of the Lie algebra $g$   
(and therefore of $\FF^{-1}$) on $X$ and $Y$
we can compose this map with the action of the twist. In this way
we obtain a deformed version $\mu_\st$ of the initial bilinear map $\mu$:
\eqa
\mu_\st:=\mu\circ \FF^{-1}~,\label{generalpres}&~~~~~~~~~~~~~~&
\ena
{\vskip -.8cm}
\eqa
{}~~~~~~~~~~~~~\mu_\st\,:X\times  Y&\rightarrow& Z\nn\\
(\x, \y)\,\, &\mapsto& \mu_\st(\x,\y)=\mu(\of^\al(\x),\of_\al(\y))\nn~.
\ena
The action of $g$ on the vectorspaces $X,Y,Z$ we consider is given by
representing $g$ as a Lie subalgebra of the Lie algebra $\Xi$ of vectorfieds
on the manifold $M$. The action of vectorfields will always be via the Lie 
derivative.
\sk
\noi 
{\bf Algebra of Functions $A_\st$}. If $X=Y=Z=Fun(M)$ where $A\equiv Fun(M)$ is the space of functions on the manifold 
$M$,
we obtain the star-product formula, for all $a,b\in A$
\eq
a\st  b=\mu\circ \FF^{-1}(a\otimes b)=\of^\al(a)\of_\al(b)~.
\en
We denote by $A_\st$ the noncommutative algebra of functions 
with the $\st$-product.
We also introduce the universal
$\RR$-matrix \eq \RR:=\FF_{21}\FF^{-1}~\label{defUR} \en where by
definition $\FF_{21}=\f_\al\otimes \f^\al$. In the sequel we use
the notation \eq
\RR=\R^\al\otimes\R_\al~~~,~~~~~~\RR^{-1}=\oR^\al\otimes\oR_\al~.
\en  The $\RR$-matrix measures
the noncommutativity of the $\star$-product. Indeed it is easy to
see that \eq\label{Rpermutation} h\st g=\oR^\al(g)\st\oR_\al(h)~.
\en The permutation group in noncommutative space is naturally
represented by $\RR$. Formula (\ref{Rpermutation}) says that the
$\st$-product is $\RR$-commutative in the sense that if we permute
(exchange) two functions using the $\RR$-matrix action then the
result does not change.
 \sk
{\bf {Basic Example. }} If $M=\mathbb{R}^n$ a main example is given by 
considering the twist
\eq
\FF=e^{-{i\over 2}\la\theta^{\mu\nu}{\partial\over \partial x^\mu}
\otimes{\partial\over \partial x^\nu}}
\en 

The $\st $-product that the twist $\FF$ induces on the algebra of
functions on $\mathbb{R}^n$ is the Moyal-Weyl $\star$-product,
\eq\label{starprod}
(f\st g)(x)=
e^{{i\over 2}\la\theta^{\mu\nu}{\partial\over \partial x^\mu}
{\partial\over \partial y^\nu}}f(x)g(y)|_{y= x}~.
\en

\sk
\noi {\bf Vectorfields $\Xi_\st$}. We now deform the $A$-module of 
vectorifelds. 
According to (\ref{generalpres}) the product 
$\mu : A\otimes \Xi\rightarrow \Xi$
is deformed into the product 
\eq
h\st v=\of^\al(h) \of_\al(v)~.
\en
Here  $\of^\al\in Ug\subset U\Xi$. Its action on vectorfiels $\of^\al(v)$  
is given by the iterated use of the Lie derivative $t(v)=[t,v]$,
 (e.g. $(t't)(v)=t'(t(v))=[t',[t,v]]$). 
It is the adjoint action. 
The cocycle condition \eqn{propF1} implies that this 
new product is compatible with the $\st$-product in $A_\st$.
We have thus constructed the $A_\st$ module of vectorfields. 
We denote it by $\Xi_\st$. 
As vectorspaces $\Xi=\Xi_\st$, but $\Xi$ is an $A$ module while $\Xi_\st$ is 
an $A_\st$ module.
\sk
\noi {\bf 1-forms $\Om_\st$}. 
The space of 1-forms $\Om$ becomes also an $A_\st$ module,
with the product between functions and 1-forms given again by following 
the general prescription (\ref{generalpres}): 
\eq
h\st\om :=\of^\al(h) \of_\al(\om)~.
\en
The action of $\of_\al$ on forms is given by iterating the Lie derivative 
action of the vectorfield $t_\mu$ on forms. 
Functions can be multiplied from the left or from the right,
if we deform the multiplication from the right we obtain the new product  
\eq
\om\st h := \of^\al(\om)\of_\al(h)
\en
and we 
``move $h$ to the right'' 
with the help of the $R$-matrix,
\eq
\omega\st h={\oR^\al}(h)\st \oR_{\al}(\om)~.
\en
We have defined the $A_\st$-bimodule of 1-forms.
  
\sk

\noi {\bf Tensorfields {$\TT_\st$}}.  Tensorfields form an algebra with 
the tensorproduct $\otimes$. We define $\TT_\st$ to be the noncommutative 
algebra of tensorfields. As vectorspaces $\TT=\TT_\st$ the noncommutative 
tensorproduct is obtained by applying (\ref{generalpres}):
\eq\label{defofthetensprodst}
\tau\otimes_\st\tau':=\of^\al(\tau)\otimes \of_\al(\tau')~.
\en 
Associativity of this product follows from the cocycle condition 
(\ref{propF1}).

\sk

\noi {\bf Exterior forms 
$\Omega^{\mbox{\boldmath $\cdot$}}_\st=\oplus_p\Omega^{p}_\st$}.
Exterior forms form an algebra with product 
$\wedge :\,\Omega^{\mbox{\boldmath $\cdot$}}\times 
\Omega^{\mbox{\boldmath $\cdot$}}\rightarrow \Omega^{\mbox{\boldmath $\cdot$}}$.
We $\st$-deform the wedge product into the $\st$-wedge product,
\eq\label{formsfromthm}
\vartheta\wedge_\st\vartheta':=\of^\al(\vartheta)\wedge \of_\al(\vartheta')~.
\en 
We denote by $\Omega^{\mbox{\boldmath $\cdot$}}_\st$ 
the linear space of forms equipped with the wedge product  
$\wedge_\st$.

As in the commutative case exterior forms are totally 
$\st$-antisymmetric contravariant tensorfields. 
For example the 2-form $\omega\wedge_\st\omega'$ is 
the $\st$-antisymmetric combination
\eq\label{stantisymm}
\omega\wedge_\st\omega'= \omega\otimes_\st\omega'
-\oR^\al(\omega')\otimes_\st \oR_{\al}(\omega)~.
\en

Since Lie derivative and exterior derivative commute, the exterior derivative 
${\rm d}:A\rightarrow \Om$ satisfies the Leibniz rule ${\rm d}(h\st g)={\rm d} h\st g+h\st{\rm d}g $, and is therefore also the $\st $-exterior derivative. 
The same argument shows that the de Rham cohomology ring is undeformed. 
\sk
\noi{\bf $\st$-Pairing}. 
We now consider the bilinear map 
$
\langle ~,~\re : \,
\Xi\times \Omega \rightarrow  A\,,$
$(v,\omega)~\mapsto \langle  v,\omega\re$,
where, using local coordinates,
$\langle  v^\mu\partial_\mu,\omega_\nu dx^\nu\re=v^\mu\om_\mu~.$
Always according to the general prescription (\ref{generalpres}) we deform 
this pairing into
\eqa\label{lerest}
\langle ~,~\re_\st : \,
\Xi_\st\times \Omega_\st &\rightarrow & A_\st~,\\
(\xi,\omega)~&\mapsto &\langle  \xi,\omega\re_\st
:=\langle \of^\al(\xi),\of_\al(\omega)\re~.
\ena
It is easy to see that due to the cocycle condition for $\FF$ 
the $\st$-pairing satisfies the  
$A_\st$-linearity properties  
\eq
\langle  h\st u,\omega\st k\re_\st=h\st\langle  u,\omega\re_\st\st k~,
\en
\eq\label{linearp}
\langle  u, h\st\omega \re_\st=
{\oR^\al}(h)\st\langle  {\oR_{\al}}(u),\omega \re_\st~.
\en
Using the pairing $\langle ~\,,~\,\res$ we associate to any $1$-form
$\om$ the left $A_\st$-linear map $\langle ~\,,\om\res$. 
Also the converse holds: any left $A_\st$-linear map 
$\Phi:\Xis\rightarrow \AAs$ is of the form $\langle ~\,,\om\res$
for some $\omega$ (explicitly $\omega=\Phi(\partial_\mu) dx^\mu$). 
\sk
The pairing can be extended to covariant tensors 
and contravariant ones. We first define in the undeformed
case the pairing $$\langle u''\ldots \otimes u'\otimes u\,,\, \theta\ldots\otimes 
\theta'\otimes\theta''\rangle =\langle u,\theta\rangle \,
\langle u', \theta'\rangle\,\langle u'', \theta''\rangle$$  ($u$, $u'$, $u''$
vectorfields, $\theta$, $\theta'$, $\theta''\,$ $1$-forms) that is obtained 
by first contracting the innermost elements. Using locality and linearity
this pairing is extended to any covariant and contravariant $n$-tensors. 
It is this onion-like structure pairing that naturally generalizes to the
noncommutative case.

The $\st$-pairing is defined by
\eq
\langle \tau, \rho \rangle_\st :=
\langle \bar\ff^\al(\tau) , \bar\ff_\al(\rho)\rangle~.
\en
Using the cocycle condition for the twist $\FF$ and the onion like structure of the undeformed pairing we have the property
\eq
\langle \tau\otimes_\st u\,,\,\theta\otimes_\st \rho \rangle_\st :=
\langle \tau \,,\,\langle u,\theta\rangle_\st\st \rho\rangle_\st~.
\en
This property uniquely characterizes the pairing
because we also have the $A_\st$ linearity properties
\eq
\langle  h\st\tau, \rho\st k\re_\st=h\st\langle  \tau,\rho\re_\st\st k~,
\en
\eq\label{linearppp}
\langle  \tau, h\st\rho \re_\st=
{\oR^\al}(h)\st\langle  {\oR_{\al}}(\tau),\rho \re_\st~.
\en
\sk
\noi {\bf $\st$-Hopf algebra of vectorfields $U\Xi_\st$}. 
\noi
Consider the universal enveloping algebra $\UU$ of vectorfields on $M$
(infinitesimal diffeomorphisms),  where $\Xi$ is the Lie algebra
of vectorfields with the usual Lie bracket 
$
[u,v](h) = u(v(h)) - v(u(h)) .
$
In order to construct the deformed algebra of vectorfields
we apply the recepy (\ref{generalpres}) and deform the product in $\UU$ 
into the new product
\eq
\xi\st\zeta=\of^\al(\xi)\of_\al(\zeta)~.
\en
We call $\UU_\st$ the new algebra with product $\st$. 
As vectorspaces $\UU=\UU_\st$.
Since any sum of
products of vectorfields in $\UU$ can be rewritten as sum of
$\st$-products of vectorfields via the formula $u\,v=\f^\al(u)\st\f_\al(v)$, 
vectorfields $u$ generate the algebra $\UU_\st$.

It turns out \cite{GR2} that $\UU_\st$ has
also a Hopf algebra structure.
We describe it by giving the coproduct, the inverse of the 
antipode and the counit on the generators $u$ of $\UU_\st$:
\eq\label{coproductu}
\D_\st (u)=u\otimes 1+ \XD_{\oR^\al}\otimes {\oR_{\al}}(u)
\en
\eq
S^{-1}_\st(u)=-\oR^\al(u)\st \XD_{\oR_\al}~.
\en
\eq
\epsi_\st(u)=\epsi(u)=0~,\label{epsist}
\en
where, for all $\xi\in\UU$,  
$\XD_\xi=\of^\al\xi\f^\be S(\f_\be) S^{-1}(\of_\al)$. The map $\XD: \UU\rightarrow\UU$ 
is invertible and it can be shown \cite{Majid-Gurevich}, 
that its inverse $\XD^{-1}$ is
\eq
\XD^{-1}_\xi=\of^\al(\xi)\of_\al=:D(\xi)~.
\en
In principle one could directly check that
(\ref{coproductu})-(\ref{epsist}) define a bona fide Hopf algebra.
Another way \cite{GR2} is to show that the Hopf algebra $\UU_\st$ is 
isomorphic to the Hopf algebra $\UU^\FF$ studied by \cite{drinfeld}. 
This latter has the same algebra structure as $\UU$. The new coproduct 
$\Delta^{\cal F}$, antipode $S^\FF$ and the counit are given by
\begin{eqnarray}
\Delta^{\cal F}(\xi) = {\cal F}\Delta(\xi){\cal F}^{-1} .\label{2.4.1bis}
\end{eqnarray}
\begin{equation}
S^\FF(\xi)=\chi S(\xi)\chi^{-1} .\label{2.4.3}
\end{equation}
\eq
\epsi^\FF(u)=\epsi(u)=0~,\label{epsistF}
\en
where $\chi := \f^\al S(\f_\al)~,\,~ \chi^{-1} = S(\of^\al) \of_\al~ .$
The isomorphism is given by the map
$D\,$: 
\eqa \label{D alg-homo}
&&D(\xi\st \zeta)=D(\xi)D(\zeta)~,\\
&&\D_\st =(D^{-1}\otimes D^{-1})\circ \D^\FF\circ D~,\label{DST}\\
&&S_\st =D^{-1}\circ S^\FF \circ D~.
\ena
\sk
Summarizing we have encountered the Hopf algebras 
$$(\UU,\cdot,\Delta,S,\epsi)~~,~~~
(\UU^\FF,\cdot,\Delta^\FF,S^\FF,\epsi)~~,~~~
(\UU_\st ,\st ,\Delta_\st ,S_\st ,\epsi)~.~$$
The first is cocommutative, the second is triangular and is 
obtained twisting the first, the third is triangular and isomorphic 
to the second. The remarkable fact about $\UU_\st$ is 
the Leibniz rule for vectorfields
(\ref{coproductu}). We have that $\oR_\al(u)$ is again a
vectorfield so that
\eq\label{woro}
\Delta_\st (\Xis) \subset \Xis \otimes 1 + \UU_\st\otimes\Xis~.
\en
This is a fundamental property for the construction of a deformed 
differential calculus \`a la Woronowicz \cite{Woronowicz}. Note that the 
coproduct $\D^\FF (u)$ does not have this property.
\sk
There is a natural action (Lie derivative) 
of $\Xi_\st$ on the space of functions $A_\st$.
It is given once again by combining the usual Lie derivative on functions 
$\ll_u(h)=u(h)$ with the twist $\FF$ as in (\ref{generalpres}),
\eq\label{stliederact}
\ll^\st_u(h):=\of^\al(u)(\of_\al(h))~.
\en
The action $\ll^\st$ of $\Xi_\st$ on $A_\st$ can be extended to all $\UU_\st$.
The map $\ll^\st$ is an action 
of $\UU_\st$ on $A_\st$, i.e. it represents the algebra $\UU_\st$ 
as differential operators on functions because
\eq
\ll^\st_u(\ll^\st_v(h))=\ll^\st_{u\st v}(h)~.\label{befmod}
\en
We also have that $\A_\st$ is a $\UU_\st$ module algebra because 
$\ll_u^\st$ is a deformed derivation of the algebra $\A_\st$ (cf. \eqn{compprim}).
Indeed in accordance with the coproduct formula (\ref{coproductu}) the differential operator $\ll^\st_u$ satisfies the deformed  Leibniz rule 
\eq
\ll^\st_u(h\st g)=\ll_u^\st(h)\st g + 
\oR^\al(h)\st \ll_{\oR_\al(u)}^\st(g)~.
\en
We conclude that $A_\st$ is a $\UU_\st$ module algebra.
\sk
This construction holds in general, see \cite{GR2}: the deformed algebras of 
functions $A_\st$, of tensorfields ${\cal T}_\st$, of exterior 
forms $\Om_\st^{\mbox{\boldmath $\cdot$}}$ and of vectorfields $\UU_\st$ are all 
$\UU_\st$ module algebras.
The $\UU_\st$ action is always given by the $\st$-Lie derivative
\eq
\ll^\st_u:=\ll_{\of^\al(u)}\circ \of^\al~.
\en 
The module property reads
\eq
\ll^\st_u\circ \ll^\st_v=\ll^\st_{u\st v}\label{modgenst}
\en
the compatibility with the algebra structure is the deformed Leibniz rule that correspond to the coproduct \eqn{coproductu}.
\sk
\noi {\bf $\st$-Lie algebra of vectorfields $\Xi_\st$}. 
In the case the deformation is given by a twist we have a natural 
candidate for the Lie algebra of the Hopf algebra $\UU_\st$.
We apply the recepy (\ref{generalpres}) and deform the Lie algebra product 
$[~,~]$ into 
\begin{eqnarray}
[\quad ]_\st: \quad\quad \Xi\times\Xi &\to& \Xi \nonumber\\
(u,v) &\mapsto& [u,v]_\st:=[\of^\al(u),\of_\al(v)]~ .\label{2.1st}
\end{eqnarray}
Notice that this $\st$-Lie bracket is just the $\st$-Lie derivative,
\eq
[u,v]_\st:=[\of^\al(u),\of_\al(v)]=\ll_{\of^\al(u)}(\of_\al(v))=
\ll^\st_{u}(v)~.
\en
In $\UU_\st$ it can be realized as a deformed commutator
\eqa 
[u,v]_\st&=&[\of^\al(u),\of_\al(v)]=\of^\al(u)\of_\al(v)-\of_\al(v)\of^\al(u)
\nn\\
 &=&u\st v-\oR^\al(v)\st\oR_\al(u)~.
\ena
It is easy to see that the bracket  $[~,~]_\st $ 
has the $\st$-antisymmetry property
\eq\label{sigmaantysymme}
[u,v]_\st =-[\oR^\al(v), \oR_\al(u)]_\st~ .
\en
This can be shown as follows
$[u,v]_\st =[\of^\al(u),\of_\al(v)]=-[\of_\al(v),
\of^\al(u)]=
-[\oR^\al(v), \oR_\al(u)]_\st~. $
A $\st$-Jacoby identity can be proven as well
\eq
[u ,[v,z]_\st ]_\st =[[u,v]_\st ,z]_\st  
+ [\oR^\al(v), [\oR_\al(u),z]_\st ]_\st ~.
\en
The appearence of the $R$-matrix $\RR^{-1}=\oR^\al\otimes\oR_\al$ is not
unexpected. We have seen that $\RR^{-1}$ encodes the noncommutativity 
of the $\st$-product $h\st g=\oR^\al(g)\st\oR_\al(h)$ 
so that $h\st g$ do $\RR^{-1}\!$-commute. Then it is natural to define  
$\st$-commutators using the $\RR^{-1}$-matrix. In other words, 
the representation of the permutation group to be used 
on twisted noncommutative spaces is the one given by the $\RR^{-1}$ matrix.

Furthermore it can be shown that the braket $[u,v]_\st$ is 
the $\st$-adjoint action of $u$ on $v$, 
\eq
[u,v]_{\st}=ad^\st_u\,v=u_{1_\st}\st v\st S(u_{2_\st})~,
\en
here we have used the coproduct notation $\D_\st(u)=u_{1_\st}\otimes u_{2_\st}$. 
More in general the $\st$-adjoint action is the adjoint action in the Hopf algebra $\UUs$, it is given by, for all $\xi,\zeta\in \UUs$,
\eq
ad^\st_\xi\,\zeta:=\xi_{1_\st}\st \zeta \st S_\st(\xi_{2_\st})~,
\en
where we used the coproduct notation 
$\D_\st(\xi)=\xi_{1_\st}\otimes\xi_{2_\st}\,$.
\sk

\sk
We call $(\Xi , [~,~]_\st)$ the $\st$-Lie algebra of vectorfields because
is a linear subspace of $\UU_\st$ such that
\eqa
i)&&\Xi_\st \mbox{ generates } \UU_\st~, \\
ii)&&\D_\st(\Xi_\st)\subset \Xi_\st\otimes 1+\UU_\st\otimes \Xi_\st ~,\\
iii)&&[\Xi_\st,\Xi_\st]_{\st}\subset \Xi_\st~.
\ena
Property $ii)$ implies a minimal deformation of the Leibniz rule. 
Property $iii)$ is the closure of $\Xi_\st$ under the adjoint action.
These are the natural conditions that according to
\cite{Woronowicz}
a $\st$-Lie algebra has to satisfy, 
see also the recent review \cite{Modave}.

\sk
\noi {\it Note 1.} The construction of the Hopf algebras $U\Xi_\st$, $U\Xi^\FF$
and of the the quantum Lie algebra $(\Xi, [~,~]_\st)$ is based only on the
twist $\FF\in U\xi\otimes U\Xi$ (the Lie derivative is just the adjoint action). Given any Lie algebra $g$ and a twist    
$\FF\in Ug\otimes Ug$ we similarly have the Hopf algebras 
$Ug_\st$, $Ug^\FF$ and the quantum Lie algebra $(g, [~,~]_\st)$.

\section{Covariant Derivative,  Torsion and Curvature}\label{covderivative}
The noncommutative differential geometry set up 
in the previous section allows to develop the formalism of covariant
derivative, torsion and curvature just by following the 
usual classical formalism.

On functions the covariant derivative equals the Lie derivative. Requiring 
that this holds in the $\st$-noncommutative case as well we immediately know 
the action of the $\st$-covariant derivative on functions, and in particular 
the Leibniz rule it has to satisfy. 
More in general we define the $\st$-covariant derivative 
$\dd^\star_u$ along the vector field $u\in \Xi$
to be the linear map $\dds_u:\Xis\rightarrow\Xis$ such that
for all $u,v,z\in\Xi_\st,~ h\in A_\st$:
\eqa
&&\dd_{u+v}^{\star}z=\dd_{u}^{\star}z+\dd_{v}^{\star}z~,\\[.35cm]
&&\dd_{h\star u}^{\star}v=h\star\dd_{u}^{\star}v~,\label{ddal}\\[.35cm]
&&\dd_{u}^{\star}(h\star v)
\,=\,\mathcal{L}_u^{\star}(h)\star v+
\oR^\al(h)\st\dd^\st_{\oR_\al(u)}v\label{ddsDuhv}
\ena
This last expression is well defined because
we have used  the coproduct (\ref{coproductu}) 
that insures that
$\oR_\al(u)$ is again a vectorfield.

The covariant derivative is extended to tensorfields using the deformed 
Leibniz rule
$$\dds_u(v\otimes_\st z)= \dds_{u}(v)\ots z + \oR^\al(v)\ots 
\dds_{\oR_\al(u)}(z)\,\,.
$$
\sk
The torsion $\tr$ and the curvature $\rr$ associated to
a connection $\dd^\st$ are the linear maps  
$\tr:\Xis \times \Xis\rightarrow\Xis$, and 
$\rr^\star : \Xis\times \Xis\times\Xis\rightarrow\Xis$ defined by
\eqa
\tr(u,v)&:=&\dd_{u}^{\star}v-\dd_{\oR^{\alpha}(v)}^{\star}\oR_{\alpha}(u)
-[u,v]_{\star}~,\\[.2cm]
\rr(u,v,z)&:=&\dd_{u}^{\star}\dd_{v}^{\star}z-
\dd_{\oR^\al{(v)}}^{\star}\dd_{\oR_\al(u)}^{\star}z-\dds_{[u,v]_\st} z~,
\ena
for all $u,v,z\in\Xis$.
{}From  the $\st$-antisymmetry property of the bracket $[~,~]_\st$, 
see (\ref{sigmaantysymme}), it
easily follows that the torsion $\tr$ and the curvature $\rr$ 
have the following $\st$-antisymmetry property
\eqa\label{Tantysymm}
\tr(u,v)&=&-\tr(\oR^\al(v),\oR_\al(u))~,\nn\\[.3em]
\rr(u,v,z)&=&-\rr(\oR^\al(v),\oR_\al(u),z)~.\nn
\ena
The presence of the $R$-matrix in the definition of torsion and curvature insures 
that  $\tr$ and $\rr$  are left $A_\st$-linear maps, i.e.
$$
\tr(f\star u,v)=f\star \tr(u,v)~
~~,~~~\tr(u,f\star v)= \oR^\al (f)\st\tr(\oR_\al(u) ,v)
$$
and similarly for the curvature.
The $A_\st$-linearity of $\tr$ and $\rr$ insures that we
have a well defined  torsion tensor and  curvature  tensor.

\sk
\noi{\bf Local coordinates description}

\noi 
We denote by $\{e_i\}$ a local frame of vectorfields 
(subordinate to an open $U\subset M$)
and by $\{\theta_j\}$ the dual frame of 1-forms:
\eq
\langle  e_i\,,\,\theta^j\re_\st=\delta^j_i~.
\en
The coefficients ${\tr_{ij}}^l$ and ${\rr_{ijk}}^l$ of the torsion and curvature tensors
with respect to this local frame are uniquely 
defined by the following expressions
\eqa
&&\tr=\theta^j\otimes_\st\theta^i\st{\tr_{ij}}^l\otimes_\st e_l\label{Tijk}~,\\
&&\rr=\theta^k\otimes_\st\theta^j\otimes_\st\theta^i
\st{\rr_{ijk}}^l\otimes_\st e_l\label{Rijkl}~,
\ena
so that
$
{\tr_{ij}}^l=\langle \tr(e_i,e_j)\,,\,\theta^l\re_\st~,~~
{\rr_{ijk}}^l=\langle \rr(e_i,e_j,e_k)\,,\,\theta^l\re_\st~.
$
We also have \cite{NG}
\eqa
&&\tr={1\over 2}_{^{}}\theta^j\wedge_\st\theta^i\st{\tr_{ij}}^l\otimes_\st e_l\label{TTijk}~,\\[.1cm]
&&\rr={1\over 2}_{^{}}\theta^k\otimes_\st\theta^j\wedge_\st\theta^i
\st{\rr_{ijk}}^l\otimes_\st e_l\label{RRijkl}~.
\ena
We now define the connection coefficients ${\Ga_{ij}}^k$, 
\eq
\nabla_{e_i}e_j=\Gamma_{ij}^k\st e_k
\label{Connectioncoeffijk}
\en 
and define the connection forms $\omega_i^j$, 
the torsion forms $\TA^l$ and the curvature forms ${\OO_k}^l$ by
\eqa
{\omega_i}^j &:=&\theta^k\st{\Gamma_{ki}}^j~,\nn\\[.1cm]
\TA^l\, &:=&-{1\over 2}_{^{}}\theta^j\wedge_\st\theta^i\st{\tr_{ij}}^l\label{THi}~,\nn\\[.1cm]
{\OO_k}^l &:=& 
-{1\over 2}_{^{}}
\theta^j\wedge_\star \theta^i\star {\rr_{kij}}^l~,\nn
\ena
It can be proven \cite{NG} that the
Cartan structural equations hold 
\eqa
&&\TA^l = d\theta^l -\theta^k\wedge_\star {\omega_k}^l \label{Cartan1}~,
\\[.1cm]
&&{\OO_k}^l = d{\omega_k}^l - {\omega_k}^m\wedge_\st{\omega_m}^l
\label{Cartan2}~.
\ena
Differentiation of the Cartan structural equations gives the 
Bianchi identities 
\eqa
&& d\TA^i +\TA^j\wedge_\st{\omega_j}^i
   =\theta^j\wedge_\star {\OO_j}^i~,\\[.1cm]
&& d\OO_k^{\ l} + {\OO_k}^m
\wedge_\star {\omega_m}^l 
- {\omega_k}^m\wedge_\star {\OO_m}^l =0~.
\ena
The proof is as in the commutative case, for example
${d}\TA^i = -{\rm d}\theta^j\wedge_\star {\omega_j}^i +
\theta^j\wedge_\star {\rm d}{\omega_j}^i 
=-\TA^j\wedge_\st{\om_j}^i-\theta^k\wedge_\st{\om_k}^j
\wedge_\st{\om_j}^i+\theta^j\wedge_\st\OO_j{}^i+\theta^k\wedge_\st\om_k{}^j\wedge_\st\om_j{}^i=-\TA^j\wedge_\st\om_j{}^i+\theta^j\wedge_\st\Om_j{}^i\,$.
\sk
\sk
We conclude this section observing that along these lines one can also consider $\st$-Riemaniann geometry.
In order to define a $\st$-metric we need to define $\st$-symmetric 
elements in $\Oms\ots\Oms$ where $\Oms$ is the space of 1-forms. 
Recalling that permutations are implemented with the $R$-matrix 
we see that $\st$-symmetric elements are of the form
\eq\label{omompr}
\omega\otimes_\st\omega'
+\oR^\al(\omega')\otimes_\st \oR_{\al}(\omega)~.
\en
In particular any symmetric tensor in
$\Om\otimes\Om\,$ is also a $\st$-symmetric tensor in 
$\Oms\ots\Oms$, indeed expansion of \eqn{omompr} gives 
$\of^\al(\om)\otimes\of_\al(\om')+\of_\al(\om')\otimes\of^\al(\om)$.
Similarly for antisymmetric tensors.

As studied in \cite{GR1,GR2} it is possible to construct a torsionfree 
metric compatible connection and to consider the appropriate $\st$-contraction 
of the Riemann tensor that leads to a well defined Ricci tensor. 
One can therefore consider Einstein equations in vacuum, 
i.e. the  vanshing of the Ricci tensor, where this last is seen 
as a function of the metric tensor.

\section{Deformed Poisson geometry}

In this section, developing the results of \cite{ALV}, we study $\st$-Poisson geometry. Let's first recall the very basic structures that we later deform.
\subsection{Poisson Bracket}
A Poisson structure on a manifold $M$ is a bilinear map
\eq\{\ ,\
\}:\mathcal A \times\mathcal A \longrightarrow \mathcal A
\label{poisson}
\en
where $\mathcal A$ is the algebra of smooth functions on $M$.
It satisfies
\beqa
&& \{f,g\}=-\,\{g,f\} ~~~~~~~~~~~~~~~~~~~~~~~~~~~~~~~~~~~~~~~{\mbox{\sl antisymmetry}}\\
&&\{f,\{g, h\}\}  + \{h,\{f,g\}\} + \{g, \{h,f\}\}= 0 ~~~~\,~~{\mbox{\sl Jacobi identity}}\\
&& \{f,g h\}= \{f, g\} h + g\{f, h\}~~~~~~~~~~~~~~~~~~~~~~~~~ {\mbox{\sl Leibniz rule}}
\eeqa

The first two properties show that the Poisson bracket $\{~,~\}$ is a Lie
bracket. We have the Lie algebra $(\A,\{~,~\})$.
We can therefore consider the
universal enveloping algebra $U\A$ that is a Hopf algebra. It is the algebra 
freely generated by the functions on $M$ modulo the ideal generated by the 
equivalence relation $f\cdot  g-g\cdot f\sim\{f,g\}$. We denote by $\cdot$ the
associative product in  $U\A$, not to be confused with the product in $\A$.

The last property (Leibniz rule) shows that the map $\{f,~\} :
\mathcal A\rightarrow \mathcal A$ is a derivation of the algebra 
of functions $\mathcal A$. In other words we have an action of the Lie algebra
$(\A,\{~,~\})$ on the algebra $\A$. 
As we recalled after \eqn{compprim} the action of 
$(\A,\{~,~\})$ on $\A$ induces an action of $U\A$ on $\A$ so that $\A$
is a $U\A$-module algebra. 

We have seen that a Poisson algebra can be equivalently defined as 
an associative algebra $\A$ that has also a compatible Lie bracket  
$\{~,~\}$. The compatibility being that the associative algebra $\A$ 
is a module algebra with respect to the Hopf algebra $U\A$.
  
Since $\A$ is the algebra of smooth functions on $M$ we have 
the Lie algebra morphism  
\eqa
X : (\A,\{~,~\})&\rightarrow &\Xi\,,\\
 f~&\mapsto &X_f := \{f, ~\}~,
\ena
$X_f$ is
the Hamiltonian vectorfield associated to the
``Hamiltonian'' $f$. 
We have  $\{f, g\}=X_f(g)$.    
In this language the Jacoby identity reads
\eq
X_{\{f,g\}}=[X_f,X_g]~.
\en
The morhism $X$ immediately lifts to the Hopf algebra morphism
\eq
X : U\A~\rightarrow ~U\Xi\,.\label{XUAU}
\en

Concerning Hamiltonian vecorfields, from  $ \{f, g\}=X_f(g)=\langle X_f,{\rm d} g \rangle $
and the antisymmetry property of the Poisson bracket we see 
that the vector field
$X_f$ actually depends on $f$ only through its differential ${\rm d}
f$, and we thus arrive at the Poisson bivector field $\La$ that
maps 1-forms into vectorfields according to
\eq
\langle\Lambda, {\rm d} f\rangle = X_f~.
\en
We therefore have
\eq
\langle\Lambda, {\rm d} f\otimes {\rm d} g\rangle = X_f(g)=\{f,g\}~.\label{onion}
\en
Notice that we use the pairing $\langle u\otimes v, {\rm d} f\otimes d
g\rangle =\langle v,{\rm d} f\rangle \,\langle u, {\rm d} g\rangle$ ($u$ and
$v$ vectorfields) that is obtained by first contracting the innermost elements.

\subsection{$\st$-Poisson Bracket}
We now consider a noncommutative manifold $(M, g,\FF,\rho)$ as defined in 
Section 3.1, and a Poisson structure on $M$.   
The natural compatibility condition between these two structures is to
require the homomorphism $\rho:g\rightarrow \Xi$ to lift to 
a homomorphism $\tilde\rho : g\rightarrow (\A,\{~,~\})$, so that 
$X\circ {\tilde\rho}=\rho$. 
If this conditon holds then (we omit writing $\rho$ and $\tilde\rho$ maps)
the twist $\FF\in U\Xi\otimes U\Xi$  is the image under 
the map $X$ in \eqn{XUAU} of a twist ${\mathscr{F}}\in U\A\otimes U\A$,
\eq
\FF=\f^\al\otimes\f_\al=X_{\ff^\al}\otimes X_{\ff_\al}~~,~~~
{\mathscr{F}}={\ff^\al}\otimes {\ff_\al}\in U\A\otimes U\A~.
\en
We can then twist the universal enveloping algebra 
$U\A$ of the Lie algebra $(A,\{~,~\})$. We therefore obtain the Hopf algebra
$U\A_\st$. The coproduct, the inverse of the 
antipode and the counit on functions 
are given by formulae \eqn{coproductu}-\eqn{epsist} where the generic 
vectorfield $u$ is replaced by the generic function $h$ and the twist 
is given by ${\mathscr{F}}$. Of course we can also consider the twisted 
Hopf algebra $U\A^{\mathscr{F}}$; it is defined as in 
\eqn{2.4.1bis}-\eqn{epsistF}.

The $\st$-Poisson bracket is then defined as (cf. \eqn{2.1st}),
\eq
\label{generaldefP}
\{f,g\}_{\star}:=\{\bar\ff^\al(f) , \bar\ff_\al(g)\} 
=\{\of^\al(f) , \of_\al(g)\} 
\en
The second equality is simply due to the fact that 
the action of ${\mathscr{F}}$ on functions is given 
by its image $\FF\in U\Xi\otimes U\Xi$.

In full analogy with the construction 
of the quantum Lie algebra of vectorfields $(\Xi, [~,~]_\st)$, we have 
the quantum Lie algebra of classical obsevables $(\A,\{~,~\}_\st)$.
In particular
the $\st$-Poisson bracket  is $\st$-antisymmetric and it satisfies
the $\st$-Jacobi identity and the $\st$-Leibniz rule\footnote{
In \cite{ALV} we do not require the twist $\FF$ to be the image of a twist 
${\mathscr{F}}\in U\A \otimes U\A$. In general we therefore do not have a Hopf
algebra $U\A_\st$. However we still have \eqn{stantis}-\eqn{leibtwist}
because we impose the milder compatibility condition between the twist $\FF$ 
and the Poisson tensor $\La$,
$\,
\bar\ff^\al\otimes \bar\ff_\al(\La) = 1\otimes \La~,~~
\bar\ff^\al(\La)\otimes \bar\ff_\al = \La\otimes
1~.
$}:
\beqa
\{f,g\}_{\star}&=&-\{\bar\rr^\alpha(g),\bar
\rr_\alpha(f)\}_{\star}~,\label{stantis}\\
\{f,\{g, h\}_\star\}_\star &=&  \{\{f,g\}_\star,h\}_\star +
\{\bar\rr^\alpha(g), \{\bar\rr_\alpha(f),h\}_\star\}_\star
\label{jacobitwist}~,\\
\{f,g\star h\}_\star &=&\{f,g\}_\star \star h +\bar
\rr^\alpha(g)\star \{\bar \rr_\alpha(f),h\}_\star
\label{leibtwist}~.
\eeqa
From the $\st$-Leibniz rule \eqn{leibtwist} we see that $\{f,~\}$
is a $\st$-vectorfield, it is indeed $\ll^\st_{X_f}$:
\eq
\ll^\st_{X_f}(g)=\ll_{\of^\al(X_f)}(\of_\al(g))=
\ll_{X_{\of^\al(f)}}(\of_\al(g))=\{\of^\al(f),\of_\al(g)\}
=\{f,g\}_\st~.
\en
In the second equality we used 
\eq
\of^\al(X_h)\otimes \of_\al= 
\of^\al(\langle \La, df\rangle)\otimes \of_\al= 
\langle \La, d\of^\al(f)\rangle\otimes \of_\al= 
X_{\bar\ff^\al(h)}\otimes
\of_\al~; \label{propfondP}
\en 
this  property holds because
 $\of^\al$ is a sum of products of Hamiltonian vectorfields.
We similarly have $
\of^\al(X_h)\otimes \of_\al=X_{\of^\al(h)}\otimes
\of_\al~$. It follows that $X$ is a morphism 
between the quantum Lie algebras 
$(\A,\{~,~\}_\st)$ and $\Xi_\st$:
\eq
[X_f,X_g]_\st=[\of^\al(X_f),\of_\al(X_g)]=[X_{\of^\al
(f)},X_{\of_\al(g)}]= X_{\{\of^\al
(f),\of_\al(g)\}}=X_{\{f,g\}_\st}~. \label{closurepoiss}
\en 
Because of property \eqn{closurepoiss}
and of the coproduct rule (\ref{coproductu}) 
we have that for a twist $\FF$ on $M$ with a 
compatible Poisson bracket {\sl Hamiltonian vector fields are
a $\star$-Lie subalgebra of the $\star$-Lie algebra of vectorfields}.
\sk
Similar techniques show that $X: U\A_\st\,\rightarrow \,U\Xi_\st$ is a homomorphism of Hopf algebras.
\sk
We can also consider the set of Hamiltonian vectorfields $\{X_Q\}$ that 
leave invariant a given Hamiltonian function $H$,
\eq
\ll^\st_{X_Q}H=0~~,{\rm i.e.},~~\{Q,H\}_\st=0 \label{symo}~.
\en
In particular, recalling \eqn{closurepoiss}, we have $[X_Q,X_H]_\st=0$.
If the Hamiltonian $H$ is invariant under the action of the twist, i.e.,
if, 
$
\of^\al\otimes \of_\al(H) = 1\otimes H~,~
\of^\al(H)\otimes \of_\al = H\otimes
1~
$
we have that $\{Q,H\}_\st=\{Q,H\}=\{H,Q\}_\st=0$ 
and $Q$ is a constant of motion.
Using the $\st$-Jacoby identity we have that 
the $\st$-bracket $\{Q,Q'\}_\st$ of two constants of motion is again a 
constant of motion. We conclude that the subspace of Hamiltonian vector 
fields $\{X_Q\}$ that leave invariant the Hamiltonian $H$ forms a 
$\st$-Lie subalgebra of the $\st$-Lie algebra of Hamiltonian vectorfields: 
the $\st$-symmetry algebra of constants of motion. The term constant of motion
is appropriate because the natural noncommutative definition of time evolution
$$\dot f=-\mathcal{L}_{X_H}^\star  f = -
\{H,f\}_\star$$ reduces to the  usual one $\dot f=-\{H,f\}$ for Hamiltonians
$H$ invariant under the action of the twist. 

\sk
\noi{\bf Example}\\
\noi
Let us consider the canonical bracket on 
phase space $M={\rm T^*}\real^n$ with the usual coordinates
$x^1,\ldots x^n, p_1,\ldots p_n$, (sum over $\ell=1,\ldots  n$ is assumed)
\eq
\{f,g\}:=\frac{\del f}{\del{x^\ell}}  ~ \frac{\del
g}{\del{p_\ell}}\,-\,\frac{\del f}{\del{p_\ell}} ~ \frac{\del
g}{\del{x^\ell}}  ~.
\en
Because of the onion like structure of the pairing and since
$\langle \frac{\del }{\del{x^i}} \otimes \frac{\del
}{\del{p_i}}\,,{\rm d} f \rangle= \frac{\del f }{\del{p_i}} \frac{\del
}{\del{x^i}}$, we have that the Poisson bivector field is
\eq
\Lambda = \frac{\del }{\del{p_i}}\,  \wedge \,\frac{\del
}{\del{x^i}} = \frac{\del }{\del{p_i}}  \otimes \frac{\del
}{\del{x^i}} - \frac{\del }{\del{x^i}} \otimes \frac{\del
}{\del{p_i}}  ~.
\en
Let $\A=C^\infty(M)$ be the space of smooth complex valued functions on $M$. 
Consider the Lie algebra $(\A,\{~,~\})$, its associated universal enveloping algebra  $(U\A,\cdot)$ where $\cdot$ is the product in $U\A$ (not to be confused with the product in $\A=C^\infty(M)$), and the twist (we absorb $\la$ in $\theta^{\ell s}$)
\eq
{\mathscr{F}}=\e^{-\frac{\ii}{2}\theta^{\ell s}
p_\ell\otimes
p_s}_{{\mbox{$\cdot$}}}{}\in U\A \otimes U\A~~,~~~~~~~~\ell, s =1,...n~.
\en
The Hopf algebra map $X:U\A\rightarrow U\Xi$ that in particular 
maps elements $f$ of the 
Lie algebra $(\A, \{~,~\})$, (i.e. functions), to the corrisponding  
Hamiltonian vectorfields $X_f$, maps the twist $\mathscr{F}$ into the twist 
\eq\label{twisttt} {\mathcal
F}=\e^{-\frac{\ii}{2}\theta^{\ell s}
\frac{\del}{\del{x^\ell}}\otimes
\frac{\del}{\del^{}{x^{s^{}}}}}\in U\Xi\otimes U\Xi
\en
%
A simple calculation shows that the twisted Poisson bracket can be
expressed as:
\eq\label{explformPoi}
\{f,g\}_{\star}= \frac{\del f}{\del{x^\ell}} \,\st\, \frac{\del
g}{\del{p_\ell}} \,-\, \frac{\del f}{\del{p_\ell}}  \,\st\,
\frac{\del g}{\del{x^\ell}} ~.
\en
Any translation invariant Hamiltonian
 is compatible with the twist and the associated
constants of motion form a $\st$-Lie algebra of constants of motions. 

We see that this formalism is quite well 
suited for field theory Hamiltonians that have potentials like 
$\int\!{\rm d}^3\!x \,\overline\Phi(x)\st\Phi(x)\st\overline\Phi(x)\st\Phi(x)$ and
are translation invariant under twists in phase space originating from twists 
on spacetime (actually on space)  like 
$\FF=\e^{-\frac{\ii}{2}\theta^{ij}\frac{\partial}{\partial x^i} 
\otimes\frac{\partial}{\partial x^j}}\,.$ This is the topic of the next section.

\section{Field Theory on noncommutative space}
\setcounter{equation}{0}
We here generalize the twist setting to the case of an infinite number of 
degrees of freedom and apply this formalism to study a scalar field theory 
on noncommutative spacetime. We choose the easiest example,
one scalar field with spacetime equal to $\real^{d+1}$
(Klein-Gordon field, or $\Phi^4$ theory). The relevant kinematical
features of classical field theories and their $\hbar$-quantization 
are already present in this simple example.  

\subsection{Classical Field Theory on noncommutative space}

The infinite dimensional phase space is described by
the fields $\Phi(x)$ and $\Pi(x)$ with $x\in\mathbb
R^d$ ($\real^{d+1}$ being spacetime). The algebra $\mathsf{A} $ is
an algebra of functionals, it is the algebra of functions on $N$
where in turn $N$ is the function space:
\eq
N={\mbox{Maps}}\,(\mathbb R^d\rightarrow\mathbb R^2)
\label{spacefields}~.
\en
We define the Poisson bracket between the
functionals $F,G\in\mathsf{A}$ to be
\eq
\{F,G\}=\int{\rm d}^{d} x \,\,\,\frac{\delta F}{\delta
\Phi}\frac{\delta G}{\delta \Pi}  - \frac{\delta F}{\delta
\Pi}\frac{\delta G}{\delta \Phi}
\en
The fields $\Phi(x)$ and $\Pi(x)$ for fixed $x$ can be considered
themselves a family of functionals parametrized by $x\in\real^n$,
for fixed $x$, $\Phi(x)$ is the functional (evaluation map) that 
associates to $\Phi$ and $\Pi$ the value $\Phi(x)$; similarly with $\Pi(x)$.
Their brackets are\footnote{In order to avoid considering
distributions we should work with smeared fields $
\Phi(f)=\int{\rm d}^{d} \!x\,\, f(x) \Phi(x)$ and $\Pi(g)=\int{\rm d}^{d}
\!x \, \,g(x) \Pi(x)$. The smeared version of the Poisson bracket
is then $ ~\{\Phi(f),\Pi(g)\}=\int{\rm d}^{d} x f(x) g(x)~. $}
\eq
\{\Phi(x),\Phi(y)\}=0~,~~
\{\Pi(x),\Pi(y)\}=0~,~~
\{\Phi(x),\Pi(y)\}=\delta(x-y)~.
\en

Now let space $\mathbb R^d$ become the noncommutative Moyal space.
The algebra of functions on $\mathbb R^d$  becomes 
noncommutative with noncommutativity
given by the twist 
$\FF=\e^{-\frac{\ii}{
2}\theta^{ij}\frac{\partial}{\partial x^i} 
\otimes\frac{\partial}{\partial x^j}}\,$
(the Lie algebra $g$ in $(M,g,\FF,\rho)$ in this case is that of translations).

The twist lifts to the algebra $\mathsf{A}$ of functionals 
\cite{Wessgauge} so that this latter too becomes noncommutative. 
This is achieved by lifting to $\mathsf A$ the action of 
infinitesimal translations.
Explicitly $\frac{\partial}{\partial x^i}$
is lifted to $\del^*_i$ acting on $\mathsf{A}$ as,
\eq
\del_i^* G:=-\int{\rm d}^d x \,\,\del_i\Phi(x)\frac{\delta G}{\delta\Phi(x)}
+\del_i\Pi(x)\frac{\delta G}{\delta\Pi(x)} 
\label{funder} ~.
\en
Therefore on functionals the twist is represented as
\eq
\mathcal F=\e^{-\frac{\ii}{2}\theta^{ij}\!\int\!{\rm d}^{\!d}\! x
\left(\del_i\Phi\frac{\delta}{\delta\Phi(x)}+\del_i\Pi\frac{\delta}{\delta
\Pi(x)}\right)\,\otimes\, \int\!{\rm d}^{\!d}\!
y\left(\del_j\Phi\frac{\delta}{\delta
\Phi(y)}+\del_j\Pi\frac{\delta}{\delta\Pi(y)}\right)}
\label{lifttwist}~.
\en
The associated $\st$-product is
\eq
F\st G= \bar\ff^\alpha(F)\bar\ff_\alpha(G)~.\label{Funstar} 
\en
If $x=y$ the  $\st$-product  between the functionals $\Phi(x)\st\Phi(y)=(\Phi\st\Phi)(x)$ where this latter 
$\st$-product is the usual one for the {\sl function} $\Phi$.

Like in the finite dimensional case, the $\st$-algebra of functionals 
$\mathsf{A}_\st$ is a $\UUs$ module algebra where now
$\Xi$ is the Lie algebra of vectorfields on $N$ 
(infinitesimal functional variations).

The vectorfield on functional space $\del_i^\st$ is a Hamiltonian vectorfield
with Hamiltonian functional $P_i\in \mathsf{A}$ given by
\eq
P_i=-\int {\rm d}^dy \,\Pi(y)\,\partial_i\, \Phi(y)~.
\en  
Indeed $
X_{P_i}=\int {\rm d}^dx {\delta P_i\over \delta\Pi(x)} {\delta \over \delta\Phi(x)}
- {\delta P_i\over \delta\Phi(x)} {\delta \over \delta\Pi(x)}=\partial^\st_i~.
$ 
The functional $P_i$ is just the momentum of the Klein-Gordon field as 
obtained via Noether theorem.
The momenta $P_i$ are mutually commuting and the twist \eqn{lifttwist} is 
therefore the image via $X\otimes X : U\A\otimes U\A\rightarrow 
U\Xi\otimes U\Xi$ of the twist $\mathscr{F}=
\e^{-\frac{\ii}{2}\theta^{i j}
P_i\otimes
P_j}_{{\mbox{$\cdot$}}}{}\in U\mathsf{A} \otimes U\mathsf{A}~.
$
We therefore twist the Hopf algebra $U\mathsf{A}$ in the Hopf algebra 
$U\mathsf{A}_\st$. The corresponding quantum Lie algebra $(\mathsf{A}, \{~,~\}_\st)$ is given by the deformed Poisson bracket, $\{~,~\}_\st~:~\mathsf{A}\otimes
\mathsf{A}\rightarrow \mathsf{A}$,
\eq
\{F,G\}_\star:= \{\bar\ff^\alpha(F),\bar\ff_\alpha(G)\}
\label{funpoi}~. \en 
This bracket satisfies the $\st$-antisymmetry, $\st$-Leibniz rule and
$\st$-Jacoby identity 
\eqa
\{F,G\}_{\star}&=&-\{\bar\rr^\alpha(G),\bar
\rr_\alpha(F)\}_{\star}\label{stantisF}\\
\{F,G\star H\}_\star &=&\{F,G\}_\star \star H +\bar
\rr^\alpha(G)\star \{\bar \rr_\alpha(F),H\}_\star
\label{leibtwistF}\\
\{F,\{G, H\}_\star\}_\star &=&  \{\{F,G\}_\star,H\}_\star +
\{\bar\rr^\alpha(G), \{\bar\rr_\alpha(F),H\}_\star\}_\star
\label{jacobitwistF}
\ena
In particular the $\st$-brackets among the fields are undeformed
\eqa
&&\{\Phi(x),\Pi(y)\}_\star=\{\Phi(x),\Pi(y)\}=
\delta(x-y)\label{PBfipi}~,\\[.4em]
&&\{\Phi(x),\Phi(y)\}_\star=\{\Phi(x),\Phi(y)\}=0~,\\[.4em]
&&\{\Pi(x),\Pi(y)\}_\star=\{\Pi(x),\Pi(y)\}=0~.
\ena
We prove the first relation 
\bea
\{\Phi(x),\Pi(y)\}_\star&=&
\{\bar\ff^\alpha(\Phi(x)),\bar\ff_\alpha(\Pi(y))\}\nonumber\\
&=&\{\Phi(x),\Pi(y)\}-\frac\ii 2 \theta^{ij}{\textstyle\left\{\int\!{\rm d}^d
z\, \del_i\Phi(z)\delta(x-z),\int\!{\rm d}^d
w\,\del_j\Pi(w)\delta(y-w)\right\}}+ O(\theta^2) \nonumber\\
&=&\{\Phi(x),\Pi(y)\}-\frac\ii 2
\theta^{ij}\del_{y^j}\del_{x^i}\delta(x-y) + O(\theta^2)
\nonumber\\
&=&\{\Phi(x),\Pi(y)\}~;
\label{4.18}\eea
the second term in the third line vanishes
because of symmetry, as well as higher terms in $\theta^{ij}$.

We conclude that for Moyal-Weyl
deformations the $\st$-Poisson
bracket just among coordinates is unchanged. It is however important to stress 
that this is not the case in general. For nontrivial functionals of the 
fields we have
\eq
\{F,G\}_\st\ne\{F,G\}   \label{nt}~.
\en

If we expand $\Phi$ and $\Pi$ in Fourier modes\footnote{We use the usual undeformed Fourier decomposition because indeed
are the usual exponentials that, once we also add the time
dependence part, solve the free field equation of motion on
noncommutative space $(\hbar^2\del^\mu\del_\mu+m^2)\Phi=0$. This
equation is the same as the one on commutative space because the
$\st$-product enters only the interaction terms.}:
\bea
\Phi(x)&=&\int\frac{{\rm d}^{d} k}{(2\pi)^{d}\,\sqrt{2E_k}}
\left(a(k)\,\e^{\ii kx} +a^*(k)\e^{-{\ii}kx}\right)\nonumber\\
\Pi(x)&=&\int\frac{{\rm d}^{d} k}{(2\pi)^{d}}
(-\ii\hbar)\sqrt{\frac{E_k}{2}} \left(a(k) \e^{{\ii} kx} -
a^*(k)\e^{-\ii kx}\right) \label{expphipi}
\eea
where $E_k=\sqrt{m^2+{\vec p}^{\,2}}=\sqrt{m^2+\hbar^2{\vec
k}^{\,2}}$, and $kx=\vec k \cdot \vec x =\sum_{i=1}^{d}k^ix^i$.
We obtain the $\st$ commutation relations between the $a(k)$ functionals
\eqa
a(k)\st a(k')
&=&\e^{-\frac{\ii}{2}\theta^{ij}\,k_ik'_j}a(k)a(k') ~~~~~,~~~~~~
a^*(k)\st a^*(k')
=\e^{-\frac{\ii}{2}\theta^{ij}k_ik'_j}\,a^*(k)a^*(k')~~, \nonumber\\[1em]
a^*(k)\st a(k')
&=&\e^{\frac{\ii}{2}\theta^{ij}k_ik'_j}\,a^*(k)a(k') ~~~~~\,,~~~~~~~
a(k)\st a^*(k')
=\e^{\frac{\ii}{2}\theta^{ij}k_ik'_j}\,a(k)a^*(k')~~. \nn
\ena
We finally calculate the Poisson bracket 
among the Fourier modes using the definition  (\ref{funpoi}) and the 
functional expressions of $a(k)$,
$a^*(k)$ in terms of $\Phi$ and $\Pi$. We obtain \cite{ALV}
\eq
\{a(k),a^*(k')\}_\st=-{i\over \hbar} (2\pi)^d\delta(k-k')~~,~~~
\{a(k),a(k')\}_\st=0~~,~~~~ \{a^*(k),a^*(k')\}_\st=0~~. 
\en
Although the twisted Poisson bracket is equal to the untwisted one 
for linear combinations of the Fourier modes,  
it yields a different result, involving nontrivial fases, 
as soon as we consider Poisson brackets of powers of $a$, $a^*$.  

\subsection{Quantum Field theory on noncommutative space}

We now formulate the canonical quantization of scalar fields on
noncommutative space.  Associated to the algebra $\mathsf{A}$ of 
functionals $G[\Phi,\Pi]$ there is the algebra $\widehat{\mathsf{A}}$
of functionals $\hat G[\hat\Phi,\hat\Pi]$ on operator valued fields. 
We lift the twist to $\widehat{\mathsf{A}}$ and then deform this algebra 
to $\widehat{\mathsf{A}}_\st$
by implementing once more 
the twist deformation principle (\ref{generalpres}).
We denote by $\hat{\partial}_i$ the lift to $\widehat{\mathsf{A}}$ 
of $\frac\partial{\partial x^i}$; 
for all $\hat G\in \widehat{\mathsf{A}}$, 
\eq
\hat\del_i \hat G:=- \int{\rm d}^d x \,\,\del_i\hat\Phi(x)\frac{\delta \hat G}{\delta\hat\Phi(x)}
+\del_i\hat\Pi(x)\frac{\delta \hat G}{\delta\hat\Pi(x)} 
\label{dhat}~;
\en
here $\del_i\hat\Phi(x)\frac{\delta \hat G}{\delta\hat\Phi(x)}$ stands for 
$\int{\rm d}^d\ell\,\del_i\Phi_\ell(x)\frac{\delta \hat G}{\delta\Phi_\ell(x)}$, where like in (\ref{expphipi}) we 
have expanded the operator $\hat\Phi(x)$ 
as $\int{\rm d}^d \ell \,\Phi_\ell(x)\hat{\sf a}(\ell)$ (and similarly for $\hat\Pi(x)$). 

Consequently the twist on operator valued functionals reads
$\hat{\mathcal F}=\e^{-\frac{\ii}{2}\theta^{ij}\hat\partial_i\otimes\hat\partial_j}$.
The twist can further be lifted to the quantum algebra of observables. 
The momentum operators are
\eq
\hat P_i=-\int {\rm d}^dy\, \mbox{\boldmath$:$}\hat\Pi(y)\,\partial_i\, \hat\Phi(y)
\mbox{\boldmath$:$}~,
\en  
where the columns $\,\mbox{\boldmath$:$}$  $\mbox{\boldmath$\;:$}\,$ stem for normal ordering.
The twist reads
$\widehat{\mathscr{F}}=
\e^{-\frac{\ii}{2}\theta^{i j}
\hat P_i\otimes
\hat P_j}\,\,.
$ 
A similar twist has been independently considered in \cite{Castro:2008su}.
With this twist $\widehat{\mathscr{F}}$ we twist the Hopf algebra $U\widehat{\mathsf{A}}$ in the Hopf algebra 
$U\widehat{\mathsf{A}}_\st$. The corresponding quantum Lie algebra 
$(\widehat{\mathsf{A}}, [~,~]_\st)$ is given by the bracket, $[~,~]_\st~:~\mathsf{A}\otimes
\mathsf{A}\rightarrow \mathsf{A}$,
\eq
[\hat F,\hat G]_\star:=[\bar\ff^\alpha(\hat F),\bar\ff_\alpha(\hat G)]
\label{funpoip}~. \en 
This deformed bracket satisfies the $\st$-antisymmetry, $\st$-Leibniz rule and
$\st$-Jacoby identity \eqn{stantisF}, \eqn{leibtwistF}, \eqn{jacobitwistF}
where all arguments are now operators and the bracket $\{~,~\}_\st$ is replaced with $[~,~]_\st$. This latter  is the $\st$-commutator. Indeed recalling the definition of the  $\RR$-matrix it can be easily 
verified that  
\eq\label{stconast}
[\hat F,\hat G]_\st= \hat F\st \hat G - 
\bar R^\alpha(\hat G)\st \bar R_\alpha(\hat F)
\en
which is indeed the $\st$-commutator in  $\hat {\mathsf{A}}_\st$.
This $\st$-commutator has been considered in \cite{Zahn}. 
\sk
We studied four Hopf algebras $U{\mathsf{A}}$, $U{\mathsf{A}_\st}$,
$U\widehat{\mathsf{A}}$, $U\widehat{\mathsf{A}}_\st$, and their corresponding
Lie  algebras
 $({\mathsf{A}},\,\{~,~\})\,,~
(\widehat{\mathsf{A}},\,[~,~])\,,~ 
({\mathsf{A}_\st},\,\{~,~\}_\st)\,,~
(\widehat{\mathsf{A}}_\st , \,[~,~]_\st)~.$ Canonical quantization on 
noncommutative space is the map $\hbar_\st$ in the diagram

\eq
\begin{CD}
{\mathsf{A}}                   @>\hbar\,\,>>  {\widehat{\mathsf{A}}}\\
 @V {\mathscr{F}} VV                                      @V \widehat{{\mathscr{F}}} VV\\
{~\mathsf{A}_\st}@>\hbar_\st >>{~\widehat{\mathsf{A}}_\st}
\end{CD}
\en
\sk
\noi We define canonical quantization 
on nocommutative space by requiring 
this diagram to be commutative as a diagram among the vectorspaces 
$\mathsf{A},\widehat{\mathsf{A}}, \mathsf{A}_\st, \widehat{\mathsf{A}}_\st$. Notice that the vertical maps, that with abuse of notation we have 
called ${\mathscr{F}}$ and $\widehat{{\mathscr{F}}}$, are the identity map, indeed $\mathsf{A}=\mathsf{A}_\st$ and $\widehat{\mathsf{A}}=\widehat{\mathsf{A}}_\st$ as 
vectorspaces. Therefore we have 
$\hbar_\st=\hbar$. 
The map $\hbar_\st$ satisfies a $\st$-correspondence 
principle because 
$\st$-Poisson brackets go into $\st$-commutators
at leading order in $\hbar$ 
\[
~~~~~~
\begin{CD}
{\{F,G\}~}                   @>\hbar~\,>>  {{-\frac{\ii}{\hbar}[\hat F,\hat G]~~}}\\[.3em]
 @V {\mathscr{F}} VV                                      @V \widehat{{\mathscr{F}}} VV\\[.2em]
\end{CD}
\]
\eq\label{diag}
\boldsymbol{
\begin{CD} ~~~~~{\{F,G\}_\st\!\!}\ @>\pmb{\hbar}_\st >>
{{-\frac{\ii}{\pmb{\hbar}}[\hat F,\hat G]_\st}}
\end{CD}}
\en
\sk
\noi Indeed recall the definitions of the $\st$-Poisson bracket and of the 
$\st$-commutator and compute
\eq
\{F,G\}_\st=\{\bar\ff^{\alpha}(F),\bar\ff_{\alpha}(G)\}\,\stackrel{\hbar}{\longrightarrow}\,-\frac{\ii}{\hbar} 
[\,\widehat{\bar\ff^{\alpha}( F)}\,,\,\widehat{\bar\ff_{\alpha}( G)}\,]
=-\frac{\ii}{\hbar} [\bar\ff^{\alpha}(\hat F),\bar\ff_{\alpha}(\hat G)]
=-\frac{\ii}{\hbar} [\hat F,\hat G]_{\st_{}}\label{5.9}
\en
The second equality holds because the lifts 
$P_i$ and $\hat{P}_i$ of $\frac\del{\del x^i}$ satisfy
$
\widehat{\{P_i,G\}}=-{\ii\over \hbar}[\widehat{P}_i,\widehat G]
$ at leading order in $\hbar$.

Repeating the passages of \eqn{4.18} we obtain  (in accordance with \eqn{5.9}) the $\st$-commutator of the fields $\hat \Phi$ and $\hat \Pi$,
\eq
[\hat \Phi(x), \hat\Pi(y)]_\st = i\hbar\delta(x-y)~. \label{fipicom}
\en
Concerning the creation and annihilation operators, they are functionals 
of the operators $\hat\Phi$, $\hat\Pi$ through the quantum analogue 
of the classical functional relation \eqn{expphipi}. 
Their $\st$-commutator follows from \eqn{fipicom} 
\eq
[\hat a(k),\hat a^\dagger(k')]_\star=(2\pi)^d\delta(k-k')~.\label{ppp}
\en

We have derived this relation from our quantization scheme for field theories on noncommutative space dictated by the symmetry Hopf  algebras $U\mathsf{A}_\st$ and $U\hat{\mathsf{A}}_\st$  and their deformed Lie brackets. In order to compare this expression with similar ones which have been found
in the literature \cite{FioreSchupp,bal1,bal2,LVV,Tureanu:2006pb,Bu,kulish,fiorewess} 
it is useful to recall \eqn{stconast} and realize the action of the 
$\RR$-matrix. Since $\mathcal R=\mathcal F^{-2}$ we obtain that (\ref{ppp})
is equivalent to
\eq
\hat a(k)\star \hat
a^\dagger(k')- \e^{-\ii\theta^{ij}k'_i k_j} \hat
a^\dagger(k')\star \hat a(k)= (2\pi)^d \delta(k-k') \label{acomm}~.
\en
This relation first appeared in \cite{Grosse}. In the noncommutative QFT 
context it appears in \cite{kulish}, \cite{Bu}, and implicitly in \cite{Zahn}, 
it is also contemplated in \cite{fiorewess} as a second option. On the other 
hand
\cite{bal1,bal2,LVV,fiorewess}, starting from a different definition of 
$\st$-commutator, $[A \stackrel{\star}{,} B]:= A\star B- B\star A$, obtain 
deformed commutation relations 
of the kind $a_k a_{k'}^\dag - \e^{-\frac{\ii}{2} \theta^{ij}
k_ik'_j} a_{k'}^\dag a_{k} = (2\pi)^d\delta(k-k') $. These are different from 
(\ref{acomm}), indeed if we expand also 
the $\st$-product in (\ref{acomm}) we obtain the usual commutation relations
$\hat a(k)\hat
a^\dagger(k')-  \hat
a^\dagger(k')\hat a(k)= (2\pi)^d \delta(k-k')$.

As in the case of the $\st$-Poisson bracket, we have found that 
the $\st$-commutator of coordinate fields \eqn{fipicom}, and 
of creation and annihilation operators \eqn{ppp}, are equal to 
the usual undeformed ones.
Once again, we warn the reader that this is not true anymore for 
more complicated functionals of the coordinate fields,  in general
$[\hat F,\hat G]_\st\not=[\hat F, \hat G]$.

\subsubsection*{Acknowledgments}
I felt very honored to present this work at the conference ``Geometry and Operators Theory'' in celebration of the 65th 
birthday of Nicola Teleman. 
This research is based on joint work with 
Marija Dimitrievi\'c, Frank Meyer, Julius Wess and with Fedele Lizzi, Patrizia Vitale;   
I would like to thank them for the fruitful collaboration.
Partial support from the European Community's Human 
Potential Program under contract MRTN-CT-2004-005104 and the
Italian MIUR under contract PRIN-2005023102 is acknowledged.

\end{document}